\documentclass[12pt]{article}
\usepackage{hyperref}
\usepackage{amsmath}
\usepackage{amssymb}
\usepackage{amsthm}
\usepackage{accents} 
\usepackage{bbm}
\usepackage{color}
\usepackage{algorithm}
\usepackage{algpseudocode}

\DeclareMathAlphabet{\mathcal}{OMS}{zplm}{m}{n}
\SetMathAlphabet{\mathcal}{bold}{OMS}{zplm}{b}{n}

\usepackage[siunitx]{circuitikz}
\usepackage{tikz}
\usepackage{diagbox}
\usepackage{pgfplots}
\usepackage{float}

\usepackage{tikz}
\usetikzlibrary{calc,patterns,decorations.pathmorphing,decorations.markings,arrows}
\usetikzlibrary{quotes,arrows.meta}

\tikzset{
  annotated cuboid/.pic={
    \tikzset{%
      every edge quotes/.append style={midway, auto},
      /cuboid/.cd,
      #1
    }
    \draw [every edge/.append style={pic actions, densely dashed, opacity=.5}, pic actions]
    (0,0,0) coordinate (o) -- ++(-\cubescale*\cubex,0,0) coordinate (a) -- ++(0,-\cubescale*\cubey,0) coordinate (b) edge coordinate [pos=1] (g) ++(0,0,-\cubescale*\cubez)  -- ++(\cubescale*\cubex,0,0) coordinate (c) -- cycle
    (o) -- ++(0,0,-\cubescale*\cubez) coordinate (d) -- ++(0,-\cubescale*\cubey,0) coordinate (e) edge (g) -- (c) -- cycle
    (o) -- (a) -- ++(0,0,-\cubescale*\cubez) coordinate (f) edge (g) -- (d) -- cycle;
    \path [every edge/.append style={pic actions, |-|}]
    (b) +(0,-5pt) coordinate (b1) 
    (b) +(-5pt,0) coordinate (b2) 
    (c) +(3.5pt,-3.5pt) coordinate (c2) 
    ;
  },
  /cuboid/.search also={/tikz},
  /cuboid/.cd,
  width/.store in=\cubex,
  height/.store in=\cubey,
  depth/.store in=\cubez,
  units/.store in=\cubeunits,
  scale/.store in=\cubescale,
  width=10,
  height=10,
  depth=10,
  units=cm,
  scale=.1,
}

\newcommand{\eg}{\emph{e.g.}}
\newcommand{\ie}{\emph{i.e.}}
\newcommand{\cf}{cf.}

\newcommand{\esssup}{\text{ess sup}}

\renewcommand{\leq}{\leqslant}
\renewcommand{\geq}{\geqslant}
\renewcommand{\epsilon}{\varepsilon}

\newtheorem{thm}{Theorem}
\newtheorem{cor}{Corollary}

\let\oldbibliography\thebibliography
\renewcommand{\thebibliography}[1]{\oldbibliography{#1}
\setlength{\itemsep}{-3pt}}

\title{\vspace{-1.3cm} A Simple Proof of Nehari's Theorem Based on Duality}
\author{Cristian R. Rojas}
\date{}

\begin{document}
\maketitle

\begin{abstract}
In this technical note we provide a simple proof of Nehari's theorem on the optimal approximation by $\mathcal{H}_\infty$ functions, based on convex duality.
\end{abstract}


\section{Introduction} \label{sec:introduction}

Nehari's theorem~\cite{Nehari-57}, relating an $\mathcal{H}_\infty$ norm minimization problem to the norm of an associated Hankel operator, is one of the cornerstones of $\mathcal{H}_\infty$ control theory. There are several proofs of this theorem in the literature: Nehari's original proof~\cite{Nehari-57,Partington-88} is based on relating the Hankel operator to a linear functional, which is then extended via the Hahn-Banach theorem to yield the result of the theorem. A different approach consists in extending the Hankel operator to a multiplication operator while preserving its norm, using Parrot's theorem; \cf~\cite{Power-82,Young1988}. Other approaches involve the commutant lifting theorem~\cite{Partington-04} and fixed point theory~\cite{Treil-Volberg-94}, among others.

Since Nehari's theorem involves a norm minimization problem, it is natural to ask whether convex duality theory can be used to establish it. This is suggested, \eg, in \cite[Sec.~10.9]{Dahleh-DiazBobillo-95}, where the authors state that ``a maximization problem in the dual space of $\mathcal{H}_\infty$ can be formulated, from which properties of the optimal solution are derived''; however, they explicitly choose not to follow that path.

In this technical note we provide a proof of Nehari's theorem based on ideas from convex duality, as inspired by D. G. Luenberger's book~\cite{Luenberger-69}. Such a proof is of great pedagogical value, as it is a natural way to approach the result. Furthermore, to the best of the author's knowledge, such a proof has not been published before, thus filling a present gap in the literature.

We assume that the reader is familiar with the notions of functional analysis at the level of \cite{Luenberger-69}. For the sake of completeness, however, we present the proofs of several preliminary results, suitably extended to complex normed spaces, because the standard required theorems, as
provided in textbooks such as \cite{Luenberger-69}, are commonly provided only for real vector spaces.


\section{Preliminaries}

In this section we present several preliminary notions and results needed to define and establish Nehari's theorem.


\subsection{Lebesgue and Hardy spaces}

Let $\mathbb{E} := \{ z \in \mathbb{C}\colon |z| > 1\}$ and $\mathbb{T} := \partial \mathbb{E} = \{ z \in \mathbb{C}\colon |z| = 1 \}$.
For $1 \leq p \leq \infty$, we define the Lebesgue space $L_p(\mathbb{T})$ as the normed space of Lebesgue measurable functions $f\colon \mathbb{T} \to \mathbb{C}$ such that\footnote{The essential supremum $\esssup f$ of a function $f\colon X \to \mathbb{R}$ is the smallest real number $a$ such that $\{x \in X\colon f(x) > a\}$ has Lebesgue measure zero.}
\begin{align*}
\begin{cases}
\displaystyle \int_{-\pi}^\pi |f(e^{i \omega})|^p d\omega < \infty, & \text{if } p < \infty \\
\displaystyle \underset{\omega \in [-\pi, \pi)}{\esssup} |f(e^{i \omega})| < \infty, & \text{if } p = \infty.
\end{cases}
\end{align*}
The norm for $L_p(\mathbb{T})$ is\footnote{We have introduced the factor $1 / (2 \pi)$ in the definitions of the $L_p$ norms to facilitate the connection between these spaces and the Hardy spaces $\mathcal{H}_p$.}
\begin{align*}
\| f \|_p :=
\begin{cases}
\displaystyle \left( \frac{1}{2 \pi} \int_{-\pi}^\pi |f(e^{i \omega})|^p d\omega \right)^{1/p}, & p < \infty \\
\displaystyle \underset{\omega \in [-\pi, \pi)}{\esssup} |f(e^{i \omega})|, & p = \infty,
\end{cases}
\end{align*}
where $f \in L_p(\mathbb{T})$. It can be shown~\cite[pp.~108]{Luenberger-69} that, for $p < \infty$ the dual of $L_p(\mathbb{T})$, \ie, the normed space of bounded linear functionals on $L_p(\mathbb{T})$, is isometrically isomorphic to $L_q(\mathbb{T})$, where $1/p + 1/q = 1$. Also, note that if $1 \leq p_1 \leq p_2 \leq \infty$, then $L_{p_1}(\mathbb{T}) \supseteq L_{p_2}(\mathbb{T})$.

The Hardy spaces~\cite{Duren-70,Hoffman-62} $\mathcal{H}_p = \mathcal{H}_p(\mathbb{E})$, for $1 \leq p \leq \infty$, are defined as the normed spaces of analytic functions\footnote{Note that we are following the control convention, by defining the Hardy spaces on $\mathbb{E}$ rather than on the open unit disk $\mathbb{D} := \{  z \in \mathbb{C}\colon |z| < 1 \}$.} $f\colon \mathbb{E} \to \mathbb{C}$ such that, for $p < \infty$,
\begin{align*}
\sup_{1 < r \leq \infty}\, \int_{-\pi}^\pi |f(r e^{i \omega})|^p d\omega < \infty,
\end{align*}
whereas, for $p = \infty$,
\begin{align*}
\sup_{z \in \mathbb{E}}\, |f(z)| < \infty.
\end{align*}
The norm for $\mathcal{H}_p$ is defined for $f \in \mathcal{H}_p$ as
\begin{align*}
\| f \|_p := \begin{cases}
\displaystyle \sup_{1 < r \leq \infty} \left( \frac{1}{2 \pi} \int_{-\pi}^\pi |f(r e^{i \omega})|^p d\omega \right)^{1/p}, & \hspace{-2mm} p < \infty \\
\displaystyle \underset{\substack{-\pi \leq \omega < \pi \\ 1 < r \leq \infty}}{\esssup}\, |f(r e^{i \omega})|, & \hspace{-2mm} p = \infty.
\end{cases}
\end{align*}
It is possible to identify $\mathcal{H}_p$, for $1 \leq p \leq \infty$, with a closed linear subspace of $L_p(\mathbb{T})$, since for every $f \in \mathcal{H}_p$ the \emph{radial limit} $\tilde{f}\colon \mathbb{T} \to \mathbb{C}$, defined as $\tilde{f}(e^{i\omega}) := \lim_{r \to 1_+} f(r e^{i \omega})$, exists for almost every $\omega \in [-\pi, \pi)$, belongs to $L_p(\mathbb{T})$ and satisfies $\|\tilde{f} \|_{L_p} = \| f \|_{\mathcal{H}_p}$; \cf~\cite[Theorem~17.11]{Rudin87} for $p < \infty$ and \cite[pp.~33]{Hoffman-62} for $p=\infty$. In particular, if $f \in \mathcal{H}_p$ for some $p$, then it can be written in the form $f(z) = a_0 + a_1 z^{-1} + a_2 z^{-2} + \cdots$ for $z \in \mathbb{E}$, and $\tilde{f}(e^{i \omega}) = a_0 + a_1 e^{-i \omega} + a_2 e^{-2 i \omega} + \cdots$ for $\omega \in [-\pi, \pi)$, so $\mathcal{H}_p$ corresponds to the closed linear subspace of $L_p(\mathbb{T})$ consisting of those $\tilde{f} \in L_p(\mathbb{T})$ such that $\int_{-\pi}^\pi \tilde{f}(e^{i \omega}) e^{-i n \omega} d\omega = 0$ for every $n \in \mathbb{N}$.

Thanks to these identifications, we will consider the $\mathcal{H}_p$ spaces interchangeably as sets of analytic functions on $\mathbb{E}$ or as subsets of $L_p(\mathbb{T})$, despite some abuse of notation.


The space $L_2(\mathbb{T})$ (and, by the identification, also $\mathcal{H}_2$) is a Hilbert space, with inner product
\begin{align*}
(f, g) := \frac{1}{2 \pi} \int_{-\pi}^\pi f(e^{i \omega})\, \overline{g(e^{i \omega})}\, d\omega, \quad f, g \in L_2(\mathbb{T}).
\end{align*}
Since $\mathcal{H}_2$ can be considered as a subspace of $L_2(\mathbb{T})$, we can define its orthogonal complement $\mathcal{H}_2^\perp$ as the linear subspace of $L_2(\mathbb{T})$ consisting of functions of the form $f(z) = a_1 z + a_2 z^2 + \cdots$, where in fact $\sum_{k=1}^\infty |a_k|^2 < \infty$. An important property of $\mathcal{H}_2$ is that if $f \in \mathcal{H}_2$ is not identically zero, then $f(e^{i \omega}) \neq 0$ for almost all $\omega \in [-\pi, \pi)$; \cf, \cite[Lemma~3.9]{Partington-88}.

By analogy, we can also define $\mathcal{H}_\infty^\perp$ as the linear subspace of $L_\infty(\mathbb{T})$ consisting of functions of the form $f(z) = a_1 z + a_2 z^2 + \cdots$. Note, however, that this is \emph{not} the orthogonal complement of $\mathcal{H}_\infty$, because there is no inner product in $L_\infty(\mathbb{T})$, and hence no notion of orthogonality in it.

Let $\mathcal{H}_1^0$ be the linear subspace of $\mathcal{H}_1$ consisting of those $f \in \mathcal{H}_1$ such that $\lim_{z \to \infty} f(z) = 0$; this means that such functions are of the form $f(z) = a_1 z^{-1} + a_2 z^{-2} + \cdots$ for $z \in \mathbb{E}$.

An important result we will need is the following:

\begin{thm}[{Riesz factorization theorem, \cite[Thm.~2.13]{Partington-88}}] \label{thm:Riesz_fact}
A function $f\colon \mathbb{E} \to \mathbb{C}$ belongs to $\mathcal{H}_1$ iff there exist functions $g,h \in \mathcal{H}_2$ such that $f = g h$ and $\| f \|_1 = \| g \|_2 \| h \|_2$.
\end{thm}

Note that if $g \in \mathcal{H}_\infty$ and $f \in \mathcal{H}_2$, then $g f \in \mathcal{H}_2$, with norm $\| g f \|_2 \leq \| g \|_\infty \| f \|_2$.


\subsection{Hahn-Banach theorem and norm minimization problems}

The Hahn-Banach theorem is a fundamental result in functional analysis. Its most popular versions apply to real normed spaces. For our purposes, we need the following extension to complex normed spaces, as presented in \cite[Ch.~3, Theorem~8]{Lax-02}:

\begin{thm} \label{thm:complex_HB}
Let $X$ be a complex vector space, and $p$ a real-valued function such that, for all $a \in \mathbb{C}$ and $x,y \in X$, $p(a x) = |a| p(x)$ and $p(x + y) \leq p(x) + p(y)$. Let $Y$ be a complex linear subspace of $X$, and $\ell$ a linear functional on $Y$ such that $|\ell(y)| \leq p(y)$ for all $y \in Y$. Then, $\ell$ can be extended to $X$ so that  $|\ell(x)| \leq p(x)$ for all $x \in X$.
\end{thm}

\begin{proof}(From~{\cite[pp.~27]{Lax-02}})
Decompose $\ell$ into its real and imaginary parts as $\ell(y) = \ell_R(y) + i \ell_I(y)$. Both $\ell_R$ and $\ell_I$ are linear over $\mathbb{R}$, and by replacing $y$ with $i y$ we see that they are related by $\ell_R(i y) = -\ell_I(y)$. Conversely, if $\ell_R$ is a real linear functional, then $\ell(x) = \ell_R(x) - i \ell_R(i x)$ is linear over $\mathbb{C}$.

Now, given $\ell$, it follows that $\ell_R(y) \leq p(y)$ for all $y \in Y$, so by the real Hahn-Banach theorem (see, \eg, \cite[Ch.~3, Theorem~1]{Lax-02}), $\ell_R$ can be extended to, say, $\ell'_R$ in $X$, so that $\ell'_R(x) \leq p(x)$ for all $x \in X$. Then, let $\ell'(x) := \ell'_R(x) - i \ell'_R(i x)$ over $x \in X$. This extended functional is linear over $\mathbb{C}$. To see that $|\ell'(x)| \leq p(x)$ holds over $x \in X$, let us write $\ell'(x) = \alpha r$, where $r \geq 0$ and $|\alpha| = 1$. Then,
\begin{align*}
|\ell'(x)|
&= r \\
&= \alpha^{-1} \ell'(x) \\
&= \ell'(\alpha^{-1} x) \\
&= \ell'_R(\alpha^{-1} x) \\
&\leq p(\alpha^{-1} x) \\
&= p(x).
\end{align*}
This concludes the proof.
\end{proof}

In particular, we will need the following corollary of the complex Hahn-Banach theorem. To introduce it, recall that a linear functional $\ell$ on a normed space $X$ is \emph{bounded} if there is an $M > 0$ such that $|\ell(x)| \leq M \| x \|$ for all $x \in X$, and the smallest such $M$ is the \emph{norm} of $\ell$, denoted as $\| \ell \|$.

\begin{cor} \label{cor:norm_preserving_extension}
Let $X$ be a complex normed space, and $\ell$ a bounded linear functional on a complex linear subspace $Y$ of $X$. Then, there exists an extension $\tilde{\ell}$ of $\ell$ to $X$ such that $\| \tilde{\ell} \| = \| \ell \|$.
\end{cor}

\begin{proof}
Take $p(x) = \| \ell \| \| x \|$, for $x \in X$, in Theorem~\ref{thm:complex_HB}.
\end{proof}

We present below a complex version of the norm minimization duality result in \cite[Sec.~5.8, Thm.~2]{Luenberger-69}. Recall that the \emph{dual space} of a normed space $X$, denoted as $X^\ast$, is the normed space of all bounded linear functionals on $X$. Also, if $x \in X$ and $x^\ast \in X^\ast$, let $\langle x, x^\ast \rangle := x^\ast(x)$.

\begin{thm} \label{thm:norm_min}
Let $M$ be a linear subspace of a complex normed space $X$. Let $x^\ast \in X^\ast$ be at a distance $d$ from $M^\perp$. Then,
\begin{align*}
d
= \min_{m^\ast \in M^\perp} \| x^\ast - m^\ast \|
= \sup_{\substack{x \in M \\ \| x \| \leq 1}} \left| \langle x, x^\ast \rangle \right|,
\end{align*}
where the minimum on the left is achieved by $m_\text{opt}^\ast \in M^\perp$. Furthermore, if $x_\text{opt} \in M$ is such that $\| x_\text{opt} \| = 1$ and $d = \langle x_\text{opt}, x^\ast \rangle$, then $\langle x_\text{opt}, x^\ast - m^\ast \rangle = \| x_\text{opt} \|\, \| x^\ast - m^\ast \|$.
\end{thm}

\begin{proof}
To show that $\min_{m^\ast \in M^\perp} \| x^\ast - m^\ast \| \geq \sup_{x \in M, \| x \| \leq 1} \left| \langle x, x^\ast \rangle \right|$, let $\epsilon > 0$ and pick an $x' \in M$ such that $\| x' \| \leq 1$ and $\left| \langle x', x^\ast \rangle \right| > \sup_{x \in M,\ \| x \| \leq 1} \left| \langle x, x^\ast \rangle \right| - \epsilon$. Then, for every $m^\ast \in M^\perp$, $\| x^\ast - m^\ast \| \geq \| x' \| \| x^\ast - m^\ast \| \geq \left| \langle x', x^\ast - m^\ast \rangle \right| = \left| \langle x', x^\ast \rangle \right| > \sup_{x \in M,\ \| x \| \leq 1} \left| \langle x, x^\ast \rangle \right| - \epsilon$, so taking $\epsilon \to 0$ proves the inequality.

To establish the reverse inequality, we need to find an $m_\text{opt}^\ast \in M^\perp$ such that
\begin{align*}
\| x^\ast - m_\text{opt}^\ast \| = \sup_{\substack{x \in M \\ \| x \| \leq 1}} \left| \langle x, x^\ast \rangle \right|.
\end{align*}
Consider the restriction of $x^\ast$ (as a linear functional) to $M$, say, $x_M^\ast$, whose norm is $\| x_M^\ast \|_M = \sup_{x \in M,\| x \| \leq 1} \left| \langle x, x^\ast \rangle \right|$. Let $y^\ast$ be a norm-preserving extension of $x_M^\ast$ to $V$, according to Corollary~\ref{cor:norm_preserving_extension}, so $\| y^\ast \| = \| x_M^\ast \|_M$ and $x^\ast - y^\ast = x_M^\ast - y^\ast = 0$ on $M$, thus $x^\ast - y^\ast \in M^\perp$. Let $m_\text{opt}^\ast = x^\ast - y^\ast$, and note that $\| x^\ast - m_\text{opt}^\ast \| = \| y^\ast \| = \sup_{x \in M,\ \| x \| \leq 1} \left| \langle x, x^\ast \rangle \right|$.

From the inequalities above, if $x_\text{opt} \in M$ is such that $\| x_\text{opt} \| = 1$ and $d = \langle x_\text{opt}, x^\ast \rangle$, then $d = \| x^\ast - m_\text{opt}^\ast \| =$ $\| x_\text{opt} \| \| x^\ast - m_\text{opt}^\ast \| \geq \left| \langle x_\text{opt}, x^\ast - m_\text{opt}^\ast \rangle \right| = \langle x_\text{opt}, x^\ast \rangle = d$, so $\langle x_\text{opt}, x^\ast \rangle = \| x_\text{opt} \| \| x^\ast - m_\text{opt}^\ast \|$, and $\langle x_\text{opt}, x^\ast - m^\ast \rangle = \| x_\text{opt} \|\, \| x^\ast - m^\ast \|$.
\end{proof}


\section{Nehari's theorem and its proof}

In this section we present Nehari's theorem, and a novel proof based on the duality result given in the previous section. Before stating the theorem, we need to introduce Hankel operators.

Let $g \in L_\infty(\mathbb{T})$. The \emph{Hankel operator with symbol $g$} is the linear operator on $\mathcal{H}_2$ given by
\begin{align*}
\Gamma_g := P_{\mathcal{H}_2} M_g R,
\end{align*}
where $P_{\mathcal{H}_2}$ is the projection of $L_2(\mathbb{T})$ onto $\mathcal{H}_2$, defined as $P_{\mathcal{H}_2}(\cdots + a_{-2} z^2 + a_{-1} z + a_0 + a_1 z^{-1} + a_2 z^{-2} + \cdots) := a_0 + a_1 z^{-1} + a_2 z^{-2} + \cdots$; $M_g$ is the \emph{multiplication operator} on $L_2(\mathbb{T})$, given by $M_g f := g f$; and $R$ is the \emph{reversion operator} on $L_2(\mathbb{T})$, defined as $R(\cdots + a_{-2} z^2 + a_{-1} z + a_0 + a_1 z^{-1} + a_2 z^{-2} + \cdots) := \cdots + a_2 z^2 + a_1 z + a_0 + a_{-1} z^{-1} + a_{-2} z^{-2} + \cdots$.

Note that, if $f_1, f_2 \in \mathcal{H}_2$ and $g \in L_\infty(\mathbb{T})$,
\begin{align*}
(\Gamma_g f_1, f_2)
&= (P_{\mathcal{H}_2} M_g R f_1, f_2) \\
&= (g [R f_1], f_2) \\
&=\frac{1}{2 \pi} \int_{-\pi}^\pi f_1(e^{-i \omega})\, g(e^{i \omega})\, \overline{f_2(e^{i \omega})}\, d\omega,
\end{align*}
where the second line follows since $M_g R f_1 = P_{\mathcal{H}_2} M_g R f_1 + (I - P_{\mathcal{H}_2}) M_g R f_1$, and $(I - P_{\mathcal{H}_2}) M_g R f_1 \in \mathcal{H}_2^\perp$, so $([I - P_{\mathcal{H}_2}] M_g R f_1, f_2) = 0$, thus $(P_{\mathcal{H}_2} M_g R f_1, f_2) = (M_g R f_1, f_2)$.

For our derivations, it will be useful to appeal to the \emph{forward shift operator} $q\colon L_1(\mathbb{T}) \to L_1(\mathbb{T})$, defined as $(q f)(z) := z f(z)$. If $f \in L_1(\mathbb{T})$, let $f^{\sim}(z) := \overline{f(z^{-1}})$; note that, for every $p \geq 1$, $f \in \mathcal{H}_p$ if and only if $f^\sim \in \mathcal{H}_p$, and $\| f^\sim \|_p = \| f \|_p$.

\medskip
We can now present our main contribution.

\begin{thm}[Nehari's Theorem] \label{thm:Nehari_theorem}
Let $g \in L_\infty(\mathbb{T})$. Then,
\begin{align*}
\inf_{f \in \mathcal{H}_\infty^\perp} \| g - f \|_\infty = \| \Gamma_g \|.
\end{align*}
Also, the infimum is attained at an $f_\text{opt} \in \mathcal{H}_\infty^\perp$. If there is a non-zero $y \in \mathcal{H}_2$ such that $\| \Gamma_g y\|_2 = \| \Gamma_g \| \, \| y \|_2$, then $g - f_\text{opt}$ is an inner function, \ie, $|g(e^{i \omega}) - f_\text{opt}(e^{i \omega})|$ is constant with respect to $\omega \in [-\pi, \pi)$.
\end{thm}

\begin{proof}
Note that $f \in \mathcal{H}_\infty^\perp$ iff $q R f \in \mathcal{H}_\infty$, because if $f(z) = a_1 z + a_2 z^2 + \cdots$ then $(q R f)(z) = a_1 + a_2 z^{-1} + a_3 z^{-2} + \cdots$, and $\esssup_{\omega \in [-\pi, \pi)} |f(e^{i \omega})| = \esssup_{\omega \in [-\pi, \pi)} |(q R f)(e^{i \omega})|$. Then,
\begin{align*}
\| g - f \|_\infty
&= \underset{\omega \in [-\pi, \pi)}{\esssup}\; |g(e^{i \omega}) - f(e^{i \omega})| \\
&= \underset{\omega \in [-\pi, \pi)}{\esssup}\; |e^{i \omega} g(e^{-i \omega}) - e^{i \omega} f(e^{-i \omega})| \\
&= \| q R g - q R f \|_\infty.
\end{align*}
Therefore, we can transform the norm minimization problem $\inf_{f \in \mathcal{H}_\infty^\perp} \| g - f \|_\infty$ into $\inf_{\hat{f} \in \mathcal{H}_\infty} \| q R g - \hat{f} \|_\infty$, where $\hat{f} = q R f$. This problem can now be put into the form of Theorem~\ref{thm:norm_min}. Indeed, by letting $X = L_1(\mathbb{T})$ (so $X^\ast = L_\infty(\mathbb{T})$), $x^\ast = q R g$, and $M = \mathcal{H}_1^0$ (so that $M^\perp = \mathcal{H}_\infty$), Theorem~\ref{thm:norm_min} establishes that the infimum is achieved by some $\hat{f}_\text{opt} \in \mathcal{H}_\infty$, and that
\begin{align*}
&\min_{\hat{f} \in \mathcal{H}_\infty} \| q R g - \hat{f} \|_\infty \\
&\qquad = \sup_{\substack{x \in M \\ \| x \| \leq 1}} \left| \langle x, q R g \rangle \right| \\
&\qquad = \sup_{\substack{x \in \mathcal{H}_1^0 \\ \| x \|_1 = 1}} \left| \frac{1}{2 \pi} \int_{-\pi}^\pi e^{i \omega}\, g(e^{-i \omega})\, x(e^{i \omega})\, d\omega \right|.
\end{align*}
In the last equality, we have replaced the inequality constraint $\| x \| \leq 1$ by the equality $\| x \|_1 = 1$, since at optimality the norm of $x$ can be shown to be as large as possible.

Note that $x \in \mathcal{H}_1^0$ iff $q x \in \mathcal{H}_1$.
Using the Riesz factorization theorem (Theorem~\ref{thm:Riesz_fact}), we can then write $q x$ as $q x = z\, w^\sim$, where $z, w^\sim \in \mathcal{H}_2$ and $\| x \|_1 = \| z \|_2\, \| w \|_2$.
Then, we have that
\begin{align} \label{eq:equalities}
&\frac{1}{2 \pi} \int_{-\pi}^\pi e^{i \omega}\, g(e^{-i \omega})\, x(e^{i \omega})\, d\omega \nonumber \\
&\qquad = \frac{1}{2 \pi} \int_{-\pi}^\pi g(e^{-i \omega})\, z(e^{i \omega})\, w^\sim (e^{i \omega})\, d\omega \nonumber \\
&\qquad = \frac{1}{2 \pi} \int_{-\pi}^\pi g(e^{i \omega})\, z(e^{-i \omega})\, \overline{w(e^{i \omega})}\, d\omega \\
&\qquad = (M_g R z, w) \nonumber \\
&\qquad = (P_{\mathcal{H}_2} M_g R z, w) \nonumber \\
&\qquad = (\Gamma_g z, w), \nonumber
\end{align}
where in the third line we performed the change of variable $\omega \leftrightarrow -\omega$, and 
the fifth line follows from the fact that $w \in \mathcal{H}_2$, so the inner product does not change its value if one projects its first argument onto $\mathcal{H}_2$. Therefore, we have that
\begin{align*}
\min_{\hat{f} \in \mathcal{H}_\infty} \| q R g - \hat{f} \|_\infty
&= \sup_{\substack{z, w \in \mathcal{H}_2 \\ \| z \|_2 = \| w \|_2 = 1}} \left| (\Gamma_g z, w) \right| \\
&= \| \Gamma_g \|.
\end{align*}
Finally, assume that there is a non-zero $y \in \mathcal{H}_2$ such that $\| \Gamma_g y \|_2 = \| \Gamma_g \| \, \| y \|_2$.
Then, letting $z = z_\text{opt} := (1/\|y\|_2)\, y$ and $w = w_\text{opt} := (1 / \| \Gamma_g y\|_2)\, \Gamma_g y$, which are both in $\mathcal{H}_2$ and have unit norm, we obtain
\begin{align} \label{eq:inner_product}
(\Gamma_g z_\text{opt}, w_\text{opt})
&= \frac{1}{\| \Gamma_g y \|_2\, \| y \|_2} (\Gamma_g y, \Gamma_g y) \nonumber \\
&= \frac{1}{\| \Gamma_g \|\, \| y \|_2^2} \| \Gamma_g y \|_2^2 \nonumber \\
&= \frac{1}{\| \Gamma_g \|\, \| y \|_2^2} \| \Gamma_g \|^2\, \|y\|_2^2 \\
&= \|\Gamma_g\|. \nonumber
\end{align}
Therefore, $(\Gamma_g z_\text{opt}, w_\text{opt}) = \sup_{z,w \in \mathcal{H}_2,\; \| z \|_2 = \| w \|_2 = 1} \left| (\Gamma_g z, w) \right|$, so $x_\text{opt}(e^{i \omega}) = e^{-i \omega} z_\text{opt}(e^{i \omega})\, \overline{w_\text{opt}(e^{-i \omega})}$ belongs to $\mathcal{H}_1^0$, has norm $1$ (as $\| x^\text{opt} \|_1 = \| z_\text{opt} \|_2 \| w_\text{opt} \|_2 = 1$), and $\langle x_\text{opt}, q R g \rangle = \sup_{x \in \mathcal{H}_1^0,\; \| x \|_1 \leq 1} \left| \langle x, q R g \rangle \right|$.
Furthermore, since $x_\text{opt}$ is the product of non-zero functions in $\mathcal{H}_2$, $x_\text{opt}(e^{i \omega}) \neq 0$ for almost all $\omega \in [-\pi, \pi)$. Then, by the last statement in Theorem~\ref{thm:norm_min}, $\langle x_\text{opt}, q R g - \hat{f}_\text{opt} \rangle = \| x_\text{opt} \|_1\, \| q R g - \hat{f}_\text{opt} \|_\infty$, or
\begin{align*}
\frac{1}{2 \pi} \int_{-\pi}^\pi x_\text{opt}(e^{i \omega})\, h(e^{i \omega})\, d\omega = \frac{1}{2 \pi} \int_{-\pi}^\pi |x_\text{opt}(e^{i \omega})|\, d\omega,
\end{align*}
where $h = (1/\| q R g - \hat{f}_\text{opt} \|_\infty) (q R g - \hat{f}_\text{opt}) \in L_\infty(\mathbb{T})$ is such that $|h(e^{i \omega})| \leq 1$ for almost all $\omega \in [-\pi, \pi)$. By the previous equation, we infer that $|h(e^{i \omega})| = 1$ for almost all $\omega \in [-\pi, \pi)$, \ie, $|e^{i \omega} g(e^{-i \omega}) - \hat{f}_\text{opt}(e^{i \omega})|$ is constant almost everywhere, which is equivalent to stating that $|g(e^{i \omega}) - f_\text{opt}(e^{i \omega})|$ is constant with respect to $\omega$ for almost all $\omega \in [-\pi, \pi)$.
This concludes the proof.
\end{proof}

The solution of the $\mathcal{H}_\infty$ norm minimization problem can be further characterized, thanks to the following corollary due to D. Sarason~\cite[Proposition~5.1]{Sarason-67} (see also \cite[Theorem~3.10]{Partington-88}):

\begin{cor}[Sarason's theorem]
Under the conditions of Theorem~\ref{thm:Nehari_theorem}, if there is a non-zero $y \in \mathcal{H}_2$ such that $\| \Gamma_g y\|_2 = \| \Gamma_g \| \, \| y \|_2$, then the minimizer $f_\text{opt}$ is unique and $g - f_\text{opt} = \Gamma_g y\, /\, R y$.
\end{cor}

\begin{proof}
Let $f_\text{opt}$ be a minimizer of $\| g - f \|_\infty$ over all $f \in \mathcal{H}_\infty^\perp$, and let $z_\text{opt}$ and $w_\text{opt}$ as in the proof of Theorem~\ref{thm:Nehari_theorem}.
Since $f_\text{opt} \in \mathcal{H}_\infty^\perp$, we have that $P_{\mathcal{H}_2} M_{g - f_\text{opt}} = P_{\mathcal{H}_2} M_g$, thus $\Gamma_{g - f_\text{opt}} = \Gamma_g$. Therefore, from \eqref{eq:inner_product}, \eqref{eq:equalities}, and the Cauchy-Schwarz inequality,
\begin{align*}
\| \Gamma_g \|
&= \| \Gamma_{g - f_\text{opt}} \| \\
&= (\Gamma_{g - f_\text{opt}} z_\text{opt}, w_\text{opt}) \\
&\leq \| P_{\mathcal{H}_2} M_{g - f_\text{opt}}\, R z_\text{opt} \|_2\, \| w_\text{opt} \|_2 \\
&\leq \| M_{g - f_\text{opt}} \|\, \| R z_\text{opt} \|_2\, \| w_\text{opt} \|_2 \\
&= \| g - f_\text{opt} \|_\infty \| z_\text{opt} \|_2\, \| w_\text{opt} \|_2 \\
&= \| \Gamma_g \|.
\end{align*}
Since the lower and upper bounds coincide, we have equality throughout, thus $M_{g - f_\text{opt}}\, R z_\text{opt} = P_{\mathcal{H}_2} M_{g - f_\text{opt}}\, R z_\text{opt} \in \mathcal{H}_2$, \ie, $\Gamma_{g - f_\text{opt}}\, z_\text{opt} =  M_{g - f_\text{opt}}\, R z_\text{opt}$, hence $g - f_\text{opt} = \Gamma_{g - f_\text{opt}}\, z_\text{opt}\, /\, R z_\text{opt} = \Gamma_g\, y\, /\, R y$ (which, in turn, implies that $f_\text{opt}$ is unique).
\end{proof}


\bibliographystyle{plain} 
\bibliography{../../../PhD/cristian.bib}

\begin{thebibliography}{10}

\bibitem{Dahleh-DiazBobillo-95}
M.~A. Dahleh and I.~J. Diaz-Bobillo.
\newblock {\em Control of Uncertain Systems: A Linear Programming Approach}.
\newblock Prentice Hall, 1995.

\bibitem{Duren-70}
P.~L. Duren.
\newblock {\em Theory of $H^p$ Spaces}.
\newblock Academic Press, 1970.

\bibitem{Hoffman-62}
K.~Hoffman.
\newblock {\em Banach Spaces of Analytic Functions}.
\newblock Prentice-Hall, 1962.

\bibitem{Lax-02}
P.~D. Lax.
\newblock {\em Functional Analysis}.
\newblock Wiley-Interscience, 2002.

\bibitem{Luenberger-69}
D.G. Luenberger.
\newblock {\em Optimization by Vector Space Methods}.
\newblock John Wiley \& Sons, 1969.

\bibitem{Nehari-57}
Z.~Nehari.
\newblock On bounded bilinear forms.
\newblock {\em Annals of Mathematics}, 65(1):153--162, 1957.

\bibitem{Partington-88}
J.~R. Partington.
\newblock {\em An Introduction to Hankel Operators}.
\newblock Cambridge University Press, 1988.

\bibitem{Partington-04}
J.~R. Partington.
\newblock {\em Linear Operators and Lnear Systems: An Analytical Approach to
  Control Theory}.
\newblock Cambridge University Press, 2004.

\bibitem{Power-82}
S.~C. Power.
\newblock {\em Hankel Operators on Hilbert Space}.
\newblock Pitman, 1982.

\bibitem{Rudin87}
W.~Rudin.
\newblock {\em Real and Complex Analysis, 3rd Edition}.
\newblock McGraw-Hill, 1987.

\bibitem{Sarason-67}
D.~Sarason.
\newblock Generalized interpolation in {$H^\infty$}.
\newblock {\em Transactions of the American Mathematical Society},
  127(2):179--203, 1967.

\bibitem{Treil-Volberg-94}
S.~Treil and A.~Volberg.
\newblock A fixed point approach to {Nehari's} problem and its applications.
\newblock In E.~L. Basor and I.~Gohberg, editors, {\em Toeplitz Operators and
  Related Topics: The Harold Widom Anniversary Volume. Workshop on Toeplitz and
  Wiener-Hopf Operators, Santa Cruz, California, September 20--22, 1992}, pages
  165--186. Springer, 1994.

\bibitem{Young1988}
N.~Young.
\newblock {\em An Introduction to Hilbert Space}.
\newblock Cambridge University Press, 1988.

\end{thebibliography}

\end{document}